\documentclass[11pt]{article}

\usepackage{latexsym}
\usepackage{amsthm}
\usepackage{amsmath}
\usepackage{amssymb}
\usepackage{amsfonts}

\newtheorem*{prop}{Proposition}

\newcommand{\la}{\langle}
\newcommand{\ra}{\rangle}
\newcommand{\bn}{\bigskip\noindent}
\newcommand{\sn}{\smallskip\noindent}
\newcommand{\bsk}{\bigskip}
\newcommand{\ssk}{\smallskip}
\newcommand{\wt}{\widetilde}
\newcommand{\bm}{\boldmath}
\newcommand{\q}{\quad}
\newcommand{\exi}{\exists \,}
\newcommand{\fa}{\forall \,}
\newcommand{\sbs}{\subseteq \,}
\newcommand{\sbsno}{\subsetneqq \,}
\newcommand{\Ry}{R[y]}
\newcommand{\My}{M[y]}
\newcommand{\kly}{[y]}

\newcommand{\ble}{\bigl[}
\newcommand{\bre}{\bigr]}
\newcommand{\bl}{\bigl(}
\newcommand{\br}{\bigr)}
\newcommand{\Z}{{\mathbb Z}}
\newcommand{\Q}{{\mathbb Q}}
\newcommand{\N}{{\mathbb N}}
\newcommand{\F}{{\mathbb F}}
\DeclareMathOperator{\Imm}{Im} \DeclareMathOperator{\imm}{im}
\DeclareMathOperator{\gker}{Ker}

\title{Computation of maximal reachability submodules}
\author{
 Wiland Schmale}
\date{
{\small Fachbereich Mathematik, Carl von Ossietzky Universit\"at,\\ D-26111 Oldenburg, Germany;\\
email: wiland.schmale@uni-oldenburg.de}}

\begin{document}
\maketitle

\sn {\bf Keywords:} Reachability submodule, linear systems, systems over rings, symbolic
computation

\begin{abstract}
A new and conceptually simple procedure is derived for the computation of
the maximal reachability submodule of a given submodule of the state space
of a linear discrete time system over a Noethenian ring $R$. The procedure
is effective if $R$ is effective and if kernels and intersections can be
computed. The procedure is compared with a rather different procedure by
Assan e.a. published recently.
\end{abstract}

\section{Introduction}

Let $ A\in R^{n\times n}$, $B\in R^{n\times m}$ where for the moment $R$ is just a commutative
ring. As usual, we associate to the pair $(A,B)$ the {\bf linear discrete time control processes}
\begin{equation}\label{eq1}
  x_0 \quad , \quad x_1 = Ax_0 + Bu_0\quad , \quad \cdots \quad ,
   \quad x_{k+1} = Ax_k + Bu_k \quad ,\quad \cdots
\end{equation}
with states $x_k\in R^n$, inputs $u_k\in R^m$ and $k\in\N$.

\bsk A submodule $U$ of $R^n$ is called {\bm$(A,B)$}{\bf -invariant}
 if $AU\sbs U + \imm B$. An
$(A,B)$ invariant submodule $U$ is called {\bf reachable} or
 {\bf reachability  submodule} if every
state in $U$ can be reached from zero within $ U$. The latter means:
\begin{multline*}
   \fa x\in U \ \exi r\in\N, \, u_0,\ldots, u_{r-1}\in R^m:\\
  x_1 = Bu_0, \ldots,x_r=A^{r-1} Bu_0 + \ldots + Bu_{r-1}\in U\ \ \text{{\bf and }} x_r=x\, .
\end{multline*}
It was shown (see e.g. \cite[Theorem 2.15]{It}) that this rather natural definition is equivalent
to the definition of pre-controllability submodules in \cite{CoPe} which is still more commonly
known but less intuitive from a control point of view.

\bsk The zero-module is trivially $(A,B)$-invariant and reachable.
 From the definitions it is clear that sums of $(A,B)$-invariant or reachable
 submodules, respectively, are
again $(A,B)$-invariant or reachable. These facts imply that any submodule $M$ of $R^n$ contains a
unique maximal $(A,B)$-invariant submodule $M^*$ and a unique maximal
reachability submodule
$M^*_0$, where always $M^*_0\sbs M^*.$

\bsk Maximal reachability submodules play an important role in the solutions to classical control
problems such as disturbance decoupling. See \cite{CoPe} and \cite{AsPe} to give only two examples.
It is therefore of practical importance to have methods at hand for the computation of generating
systems of such modules. In \cite{AsLaPe1} for the first time a finite  procedure was given for
principal ideal domains and then strongly modified in \cite{ AsLaPe2} to work for Noetherian rings.
The latter works as follows:

\bn $R$ is now supposed to be Noetherian.
\begin{description}
\item[First step (precalculation):]
   $S_0:= \imm B$ \\ and for $k\ge 1:\quad S_k := \imm B + A(S_{k-1}\cap M)$.\\
This ascending sequence of modules stabilizes after finitely many steps and gives a submodule $M_*$
which contains the image of $B$. If $M$ is represented as the kernel of some matrix $C\in
R^{n\times p}$, then $M_*$ appears as the 'minimal $(C,A)$-invariant submodule' containing the
image of $B$, see e.g. \cite{AsLaPe2}.

\item[Second step and main procedure:]
$   W_0 := M_*\cap M \cap A^{-1} (\Imm B)$ \\
 and for $k\ge 1$:\quad $W_k := M_*\cap M \cap A^{-1}
(W_{k-1} + \Imm B)$.\\ Once more, this gives an ascending sequence and an interesting proof in
\cite{AsLaPe2} shows that its limit is actually $M^*_0$.
\end{description}

\sn Of course - and the same is valid for the new procedure to be developed in this note - such a
procedure can be realized in a concrete computation only if the ring $R$ and all the occurring
operations like ``$A{^{-1}}$'', ``$\cap$'' are  effective in the sense of \cite[p.1]{CoCuSt}.

\section{New procedure via finite {\bm $(A,B)$}-cyclic\\ submodules}

Based on results from \cite{BrSch} a quite different and conceptually simpler approach is possible.
A submodule $U$ of $R^n$ is called {\bm $(A,B)$}{\bf -cyclic}  if for some $u_k\in R^m$ and $x_k$
from \eqref{eq1} with \quad {\bm $ x_o=0$} \quad one has

\begin{equation}\label{eq2}
   U = \la x_k : k\ge 0\ra\, .
\end{equation}
Thus an $(A,B)$-cyclic submodule can be generated by the states of one single control process which
begins with the zero-state.

It is shown in \cite{BrSch} that $(A,B)$-cyclic submodules are reachability submodules and that
{\bf finitely generated} reachability submodules are even {\bf finite} $(A,B)$-cyclic. The latter
means that in addition to \eqref{eq2} one has $x_k=0$ for $k>d$ and some $d\in\N$.

\bsk The point is now that finite $(A,B)$-cyclic submodules can be determined via the kernel of $[
yE-A, -B]$ in $\Ry^{n+m}$. If for $f\in \Ry^n$, $g\in \Ry^m$ one has $(yE-A)f = Bg$, then the
coefficient vectors of $f$ generate a finite $(A,B)$-cyclic submodule and every finite
$(A,B)$-cyclic submodule $U= \la x_1,\ldots, x_d, \, 0,\ldots\ra$ leads to a kernel element $\left[
\begin{smallmatrix} f\\g
\end{smallmatrix}\right] $
with $f = x_1y^{d-1}+\ldots + x_d$ and $g= u_0 y^d +\ldots + u_d$. Note that $x_{d+1} = A_dx_d+
Bu_d=0$. More details can be found in \cite{BrSch}.

For any $f= x_1y^{d-1} + \ldots + x_d\in R\kly^n$ let \quad {\bm $ U_f: = \la x_1,\ldots,
x_d\ra$}.\\ Of course, $U_f$ is contained in a given submodule $M$ if and only if the coefficient
vectors of $f$ are from $M$. Let $\pi$ be the projection of $R\kly^{n+m} = \Ry^n \oplus \Ry^m$ onto
the first $n$ components and let
\begin{equation}\label{eq3}
   {\mathcal M} := \gker \ble yE-A,-B\bre \cap \ \bl \My \times \Ry^m\br\, .
\end{equation}
Here $M[y]$ is the submodule of $R[y]^n$ of those polynomial vectors which have all their
coefficient vectors from $M$.

 One arrives now at the following results:

\bn {\bf Observation}. (i) For every $h\in {\mathcal M}$ the submodule $U_ {\pi(h)}$ is a
reachability submodule of $M$ (true for any $R$).

\ssk (ii) Let $R$ be Noetherian. For every reachability submodule $U$ of $M$ there is $h\in
{\mathcal M}$ such that $U= U_{\pi(h)}$.

\bn \begin{prop} Let $h_1,\ldots, h_s$ generate ${\mathcal M}$ as an $\Ry$-module, then the family
of coefficient vectors of $\pi(h_1),\ldots, \pi(h_s)$ generates $M^*_0$.
\end{prop}

\sn
\begin{proof}[{\bf Proof of Observation}]
 (i): By construction $U_{\pi(h)}$ is finite $(A,B)$-cyclic and thus by
Proposition 1.5 in \cite{BrSch} a reachability submodule.

\sn (ii): Since $R$ is Noetherian, $U$ is finitely generated and reachable. By Proposition 1.7 in
\cite{BrSch} this implies that $U$ is finite $(A,B)$-cyclic. The foregoing discussion shows how to
construct the desired $h\in {\mathcal M}$.
\end{proof}

\begin{proof}[{\bf Proof of Proposition}]
Let $f_1=\pi(h_1), \ldots, f_s=\pi(h_s)$ and $\wt{M} = \sum^s_{i=1} U_{f_i} $. We have to show
$\wt{M}=M^*_0$. $M^*_0$ is the sum of all reachability submodules of $M$. Since $R$ is Noetherian,
all reachability submodules $U$ of $M$ are finitely generated. By part (ii) of the Observation such
modules $U$ can be represented as $U=U_{\pi(h)}$ with some $h\in {\mathcal M}$. Since
$h=r_1h_1+\ldots + r_sh_s$ with some $r_1,\ldots, r_s\in\Ry$, we obtain $U\sbs \wt{M}$ for an
arbitrary reachability submodule $U$  of $M$ and thus $M^*_0\sbs \wt{M}$.

The converse inclusion comes from the fact that by part (i) of the Observation $U_{f_i}$ is a
reachability submodule of $M$ and therefore contained in $M^*_0$ for $1 \le i\le s$. The latter
implies: $\wt{M}\sbs M^*_0$.
\end{proof}

One main advantage of the approach via \eqref{eq3} is that one can (for appropriate rings $R$)
compute the kernel of $[yE-A,-B]$ once for all independently of $M$. This gives us  as a first
result a module which is of use not only for  determining $M^*_0$, see e.g. \cite{BrSch}. In order
to determine $M^*_0$ for some specific $M$ it remains to calculate an intersection of two modules
and after that one merely truncates the results and extracts the coefficient vectors.. Explicit
calculation is - of course - only possible over an effective Noetherian ring with an effective
method to determine the kernel and intersection in \eqref{eq3}. Examples of such rings are $\Z,\Q\,
[t_1,\ldots, t_n]$, $\F\,[t_1,t_n]$ where $\F$ is a finite field. The determination of
$\gker\,[yE-A,-B]$ can then be done with the help of Gr\"obner basis calculations as indicated in
\cite{BrSch}. A standard technique also via Gr\"obner bases for the computation of the
intersections of modules is (e.g.) described in \cite{AdLou}. In both cases any generating system
would do as well. Several current software packages for symbolic computation can be utilized to
perform explicit calculations.

 A sound comparison of the different procedures for the computation of maximal reachability
submodules requires a detailed investigation of their complexities. This remains as a future task.

\bsk The following two examples are over $\Q[t]$ and $\Q[t,w]$. Computations have been done
combining the well-known packages Macaulay2 and MapleV Release 5.1 \bsk

{\bf Examples \\

(A)}\quad  Let $A= \left[ \begin{smallmatrix} 0 & 1 & 0\\ 0 & 0 & t\\ 0 & 0 & 0
\end{smallmatrix} \right] \;
B=
\left[ \begin{smallmatrix}
1 & -t\\ t & t \\ 0 & t
\end{smallmatrix} \right] \;
$
and $M=\imm$
$\left[ \begin{smallmatrix}
-1 & 0\\ 1 & 0 \\ 0 & 1
\end{smallmatrix} \right] $
as in \mbox{Example 1 } of \cite{AsLaPe2}.

 To determine $M^*_0$ we first obtain
\[
   \gker [yE-A,-B] = \imm
     \begin{bmatrix}
         t & -t-y\\ -t & -ty\\ -t & 0\\t & -y^2\\ -y & 0
     \end{bmatrix}
\]

\noindent This leads to ${\mathcal M}=h \Ry$ with $h = \, ^t[t,-t,-t,t,-y]$, which in turn leads to
with $f=\pi(h) = \,^t[ t,-t,-t]$. There is only one coefficient vector to be extracted from $f$
(viewed as a polynomial vector in the variable y). Therefore the final result is:\; $M^*_0=fR$.\,
By \cite{AsLaPe2} we know $M^*= \left[
\begin{smallmatrix} 1
\\ -1 \\ -1
\end{smallmatrix} \right]  R$ \
and thus $M^*_0\sbsno M^*$ .

This example is interesting also since here the classical Wonham-algorithm to determine $M^*$ does
not converge and up to now no general finite procedure is known. For principal ideal domains,
however, a procedure has been developed in \cite{AsLaLoPe}.

\bn {\bf (B)}\quad  In the second example we start with matrices from \cite{AsLaPe2}, Example 4.3,
where a system with two incommensurable delays is investigated.

\bn
Let
\[
A=
\begin{bmatrix}
  0 & 0 & 1\\ w^4 & t & 0\\x^3 & t & 1
\end{bmatrix}\; ,
\; B=
\begin{bmatrix}
  t& 0 \\ 0& 0\\0& 1
\end{bmatrix}\; ,
  \text{ and }  M = \imm
\begin{bmatrix}
  1& 0 \\ 0& w\\0& 1
\end{bmatrix}\; .
\]

\sn
Here Macaulay2 computes
\[
  \gker\,[yE-A,-B] = \imm
\begin{bmatrix}
  0&\q -t+y \\ 0&\q w^4\\-t&\q -ty+y^2\\
    1&\q 0\\ -ty &\q (-w^4t+t^4)-t^3y-ty^2+y^3
\end{bmatrix}\; ,
\]
which leads to ${\mathcal M}=h\Ry$ with
\[
  h=\, ^t\ble t^2-ty,\, -w^4t, -w^3t,\, -w^3 + ty-y^2,\, (w^4t^2-t^5) +
    (-w^3t+t^4)y\bre\, .
\]
Now $\pi(h) = x_1y+x_2$ where $x_1=\, ^t[-t,0,0]$ and $^tx_2 =$ $[ t^2,-w^4t,-w^3t]$ and according
to the Proposition we obtain as final result:\quad$M^*_0=\la x_1,x_2\ra$ (compare with $R^*_2$ in
\cite[4.3]{AsLaPe2}). Note that by the new procedure we automatically get $M^*_0$ represented as an
$(A,B)$-cyclic subspace. In more complex examples one obtains $M^*_0$ as a sum of $(A,B)$-cyclic
modules. For reasons of space I do not give an example for this.

\end{document}